\documentclass[journal]{IEEEtran}

\IEEEoverridecommandlockouts     

\usepackage{color,amsmath,amstext,amsfonts,amssymb, mathrsfs,mathtools}
\usepackage{graphicx,algorithm}
\usepackage[algo2e,linesnumbered]{algorithm2e}
\usepackage{tikz}
\usepackage{setspace}
\usepackage{array} 
\usepackage{booktabs}  
\usepackage{bbm}
\usepackage{wasysym}
\usepackage{textcomp}
\usepackage{hyperref} 
\usepackage[sort,compress]{cite}

\newtheorem{lemma}{Lemma}
\newtheorem{assumption}{Assumption}

\newtheorem{theorem}{Theorem}
\newtheorem{definition}{Definition}
\newtheorem{corollary}{Corollary}

\renewcommand{\thesubsubsection}{{\it \Alph{subsection}.\arabic{subsubsection}.}}

\SetKwInput{KwInput}{Input}

\newcounter{l1}
\newcounter{l2}
\newcounter{l3}
\newcommand{\bdotlist}{\begin{list}{$\bullet$}{}}
\newcommand{\bboxlist}{\begin{list}{$\Box$}{}}
\newcommand{\bbboxlist}{\begin{list}{\raisebox{.005in}{{\tiny $\blacksquare$ \ \ }}}{}}
\newcommand{\bdashlist}{\begin{list}{$-$}{} }
\newcommand{\blist}{\begin{list}{}{} }
\newcommand{\barablist}{\begin{list}{\arabic{l1}}{\usecounter{l1}}}
\newcommand{\balphlist}{\begin{list}{(\alph{l2})}{\usecounter{l2}}}
\newcommand{\bAlphlist}{\begin{list}{\Alph{l2}.}{\usecounter{l2}}}
\newcommand{\bdiamlist}{\begin{list}{$\diamond$}{}}
\newcommand{\bromalist}{\begin{list}{(\roman{l3})}{\usecounter{l3}}}

\providecommand{\norm}[1]{\lVert#1\rVert}

\newcommand{\beq}{\begin{equation}}
\newcommand{\eeq}{\end{equation}}

\newcommand{\tn}{\textnormal}

\title{Adaptive Gradient Online Control}

\author{Deepan Muthirayan, Jianjun Yuan, Pramod P. Khargonekar
\thanks{This work is supported in part by the National Science Foundation under Grant ECCS-1839429.
D. Muthirayan and P. P. Khargonekar are with the Department of Electrical Engineering and Computer Sciences, University of California Irvine, Irvine, CA (emails: deepan.m@uci.edu, pramod.khargonekar@uci.edu). Jianjun Yuan is with the Expedia Group (email:yuanx270@umn.edu).}
}

\begin{document}

\maketitle
\thispagestyle{empty}
\pagestyle{empty}

\begin{abstract}
In this work we consider the online control of a known linear dynamic system with adversarial disturbance and adversarial controller cost. The goal in online control is to minimize the regret, defined as the difference between cumulative cost over a period $T$ and the cumulative cost for the best policy from a comparator class. For the setting we consider, we generalize the previously proposed online Disturbance Response Controller (DRC) to the adaptive gradient online Disturbance Response Controller. Using the modified controller, we present novel regret guarantees that improves the established regret guarantees for the same setting. We show that the proposed online learning controller is able to achieve intermediate intermediate regret rates between $\sqrt{T}$ and $\log{T}$ for intermediate convex conditions, while it recovers the previously established regret results for general convex controller cost and strongly convex controller cost.
\end{abstract}


\section{introduction}

Control of systems with uncertainties is a central challenge in control and is an extensively researched topic. There are various sub-fields in control such as stochastic control \cite{kumar2015stochastic, aastrom2012introduction}, robust control \cite{skogestad2007multivariable} and adaptive control \cite{sastry2011adaptive, ioannou2012robust} that address the challenge of controller synthesis for different types of uncertainties. 
In this work we are concerned with the problem of online control of systems with uncertainties such as disturbance and adversarial controller cost. The performance in online control is measured in terms of how the regret of performance, defined as the deviation of the performance from that of the best policy, scales with the duration $T$. The objective in online control is to design adaptive algorithms to disturbances and adversarial cost so that the regret scales sub-linearly in $T$, i.e., as $T^\alpha$ with $\alpha < 1$. 

Classical adaptive control investigates the problem of control of systems with parametric, structural and parametrizable disturbance uncertainties \cite{tao2014multivariable}. The main focus in classical adaptive control is the stability of system and asymptotic tracking performance. Adaptive control has been studied for systems of all types such as linear, non-linear, and stochastic. There are many variants of adaptive control such as adaptive model predictive control \cite{heirung2017dual, lorenzen2017adaptive}, adaptive learning control \cite{marino2012robust, yu2015switching}, stochastic adaptive control \cite{aastrom2013adaptive} and robust adaptive control \cite{ioannou2012robust}. These variations address the design of adaptive controller for different variations of the basic adaptive control setting. Many papers and books have been written on adaptive control; see for example  \cite{sastry2011adaptive, ioannou2012robust, aastrom2013adaptive}. Thus, adaptive control is a very rich and extensively studied topic. The key variation of the online control setting from the classical adaptive control is the regret objective and in some cases the general nature of the costs, where they could be adversarial and unknown apriori. Thus, the classical adaptive control approaches can be inadequate to analyse online control problems and are typically solved by merging tools from statistical learning, online learning and optimization, and control theory.

The field of online control has seen rising interest in the last few years. One of the first setting that was extensively explored is the Linear Quadratic Regulator (LQR) with the unknown system and stochastic disturbances. Abbasi \& Czepesvari \cite{abbasi2011regret} were the first to study the online LQR problem with unknown system and stochastic disturbances. 
The authors proposed an adaptive algorithm that achieved $\sqrt{T}$ regret w.r.t the best linear control policy, which is the optimal policy. After \cite{abbasi2011regret}, several authors improved the algorithm of \cite{abbasi2011regret}, which was an inefficient algorithm. Dean et al. \cite{dean2018regret} were the first to propose an efficient algorithm for the same problem. They showed that their algorithm achieved a regret of $\mathcal{O}(T^{2/3})$. Cohen et al. \cite{cohen2019learning} and Mania et al. \cite{mania2019certainty} improved on this result by providing an efficient algorithm with a regret guarantee of $\mathcal{O}(T^{1/2})$ for the same problem. Mania et al. \cite{mania2019certainty} extended these results to the partial observation setting and established $\mathcal{O}(\sqrt{T})$-regret for the partially observed Linear Quadratic Gaussian (LQG) setting. Cohen et al. \cite{cohen2018online} provided an $\mathcal{O}(\sqrt{T})$ algorithm for a variant of the online LQR, where the system is known and noise is stochastic but the controller cost function is an adversarially chosen quadratic function. Recently, Simchowitz et al. \cite{simchowitz2020naive} showed that $\mathcal{O}(T^{1/2})$ is the optimal regret for the online LQR control problem. 

While the above works focussed on online LQR, there are others who studied the control of much general systems: linear dynamic systems with adversarial disturbances and adversarial cost functions. Agarwal et al. \cite{agarwal2019online} considered the control of a known linear dynamic system with additive adversarial disturbance and an adversarial convex controller cost function. They proposed an online learning algorithm that learnt a Disturbance Response Controller (DRC): a linear feedback of the portion of the output contributed by the disturbances upto certain history. They showed that their proposed controller achieves $\mathcal{O}(\sqrt{T})$-regret with respect to the best DRC in hindsight. Agarwal et al. in a subsequent work \cite{agarwal2019logarithmic} showed that a poly logarithmic regret is achievable for strongly convex controller cost and well conditioned stochastic disturbances. Hazan et al. \cite{hazan2020nonstochastic} extended the setting of \cite{agarwal2019online} to the case where the system is unknown. They showed that when the system is unknown, while $\mathcal{O}(\sqrt{T})$-regret is not achievable, they can still achieve a sub-linear regret of $\mathcal{O}(T^{2/3})$-regret. Recently, \cite{simchowitz2020improper} generalized these results to provide similar regret guarantees for the same setting with partial observation for both known and unknown systems.

In this work we study the online control setting of \cite{simchowitz2020improper}: linear dynamic systems with additive disturbance and adversarial controller cost, where the system state is only partially observable. We assume that our system is known and our cost functions are general convex controller costs. Previous works in the online adversarial setting \cite{agarwal2019online, agarwal2019logarithmic, hazan2020nonstochastic, simchowitz2020improper}, either assume the cost functions to be convex or strongly-convex. Reiterating the results of \cite{simchowitz2020improper} for the known system case, what has been established is that $\mathcal{O}(\sqrt{T})$ regret is achievable when the cost function are convex, and $\mathcal{O}(\log{T})$ regret is achievable when the cost functions are strongly convex. The question we address in this work is: {\it can we achieve intermediate regret guarantees for intermediate convex conditions}?

\subsection{Our Contribution}

The online control algorithm we propose is the adaptive gradient extension of the online learning disturbance response controller proposed in \cite{agarwal2019online, simchowitz2020improper}. Here the adaptive gradient refers to the adaptation of the gradient step size of the gradient learning algorithm used in \cite{agarwal2019online, simchowitz2020improper}. Thus, to the best of our knowledge, we present the {\it first adaptive gradient online learning control algorithm.} We show that the proposed learning algorithm {\it recovers the previously established regret guarantee of $\mathcal{O}(\sqrt{T})$ for general convex controller cost functions and $\mathcal{O}(\log{T})$ for strongly-convex and smooth controller cost functions (see \cite{simchowitz2020improper}), and simultaneously achieves an intermediate regret between $\mathcal{O}(\sqrt{T})$ and $\mathcal{O}(\log{T})$ for intermediate convex conditions of the controller cost functions}. We prove our main result by establishing a new result for adaptive gradient online learning for the problem of Online Convex Optimization with Memory (OCO-M), which is the online convex optimization problem where the cost at a time step also depends on a certain history of past decisions.

\subsection{Other Related Work}

{\it Online Convex Optimization (OCO)}: In the OCO framework, the learner encounters a sequence of convex loss functions which are unknown beforehand and may vary arbitrarily over time. The learner updates the estimate of the optimal solution at each time-step based on the previous losses and incurs a loss for its updated estimate as given by the loss function for this time step. At the end of each step, either the loss function may be revealed, a scenario referred to as full information feedback, or only the experienced loss is revealed, a scenario known as bandit feedback. The objective of the learner is to minimize the loss accumulated over time. Under the full information feedback setting, it has been established that the best possible regret scales as $O(T^{1/2})$ (resp. $O(\log T)$) for convex (resp. strongly convex) loss functions, where $T$ is the number of time steps \cite{zinkevich2003online, hazan2006logarithmic, abernethy2009stochastic}. These results have also been extended to constrained online convex optimization where it has been shown that the best regret scales as $O(T^{\max\{c,1-c\}})$ for the cost and $O(T^{1-c/2})$ for constraint violation, where $c$ is a constant \cite{jenatton2016adaptive, yuan2018online}. When compared to OCO, the key difference in online control is the dependence of the decision on the state of the system, and thus in online control what is to be learnt is a control policy instead of a single decision.

{\it Policy Optimization}: Fazel et al. \cite{fazel2018global} proved that the policy gradient based learning converges asymptotically to the optimal policy for the Linear-Quadratic Regulator (LQR) problem. Zhange et al.  \cite{zhang2019policy} extended this result to the ${\cal H}_2/{\cal H}\infty$ control problem. Recently, \cite{molybog2020global} proved asymptotic convergence of a gradient based meta-learner for the LQR problem. All of these works provide asymptotic convergence guarantees.

{\it Notation}: We denote the transpose of a vector $X$ by $X^\top$. We denote the expectation of a random variable $X$ by $\mathbb{E}[X]$ and the expectation w.r.t a filtration $\mathcal{F}_t$ by $\mathbb{E}[. \vert\mathcal{F}_t]$. The minimum singular value of a matrix $M$ is denoted by $\sigma_{\text{min}}(M)$ and the minimum eigen value is denoted by $\lambda_{\text{min}}(M)$. The function $\rho(\cdot)$ denotes the spectral radius of the input matrix. We define $\norm{\cdot}$ to be 2-norm of the vector or the matrix as the case maybe. For a given variable $X_{t}$ that is dependent on time $t$, $X_{t_1:t_2}$ is used to denote the sequence $(X_{t_1}, X_{t_1+1},...,X_{t_2})$. By $\sum X_{t1:t2}$, we denote the sum of the elements in the sequence $X_{t_1:t_2}$. The big $\mathcal{O}(\cdot)$ is the standard order notation and $\tilde{\mathcal{O}}(\cdot)$ is the standard order notation that includes polylog factors. 

\section{Problem Preliminaries}

\noindent The problem we consider is the online control of a linear dynamical system given by
\begin{align}
& x_{t+1} = Ax_t +Bu_t +w_t, \nonumber \\
& y_t = Cx_t + e_t,
\label{eq:sys-dyn}
\end{align}
where $x_t \in \mathbb{R}^{d_x}$, is the state of the system, $u_t \in \mathbb{R}^{d_u}$, is the control input generated by the controller, $w_t, e_t$ are bounded disturbances of appropriate dimensions and $y_t \in \mathbb{R}^{d_y}$, is the observed output. The objective is to regulate the response of this system so as to achieve sub-linear regret with respect to the best policy from a class of policies, also called the comparator policy.

The class of policies we consider for the comparator are {\it linear dynamic controllers}, denoted by $\Pi$. A linear dynamic controller $\pi \in \Pi$ is a linear dynamic system given by $(A_\pi, B_\pi, C_\pi, D_\pi)$ with the internal state $s^\pi_t \in \mathbb{R}^{d_\pi}$ and output being the control input at time $t$:
\begin{equation}
s^\pi_{t+1} = A_\pi s^\pi_t +B_\pi y_t, u^\pi_t = C_\pi s^\pi_{t} + D_\pi y_t
\label{eq:ldc}
\end{equation}

We denote the online controller for the system in Eq. \eqref{eq:sys-dyn} by $\mathcal{C}$. The controller at any point has only access to the following information at time $t$: (i) all prior cost functions $c_{1:t-1}$, (ii) all prior observations $y_{1:t-1}$, and (iii) all prior control inputs $u_{1:t-1}$. The controller, unlike the classical setting, does not have access to the future cost functions, which are adversarial. The controller has to choose a policy to compute the control action at time $t$ based on this information.

{\it The online control setting of ours is the following}: The controller $\mathcal{C}$, on applying the control input $u_t$ at time $t$, suffers the loss $l_t(y_t,u_t)$, an adversarially chosen convex function, which is apriori unknown. The controller can observe the loss function only after its decision at time step $t$. The controller can then use this information to update its control policy. 

The performance of the online controller is measured by the regret which is the total cost incurred by the controller for a duration $T$ minus the total cost incurred by the best controller in hindsight taken from the class of controllers $\Pi$. Denote the system output and the input corresponding to a controller $\pi \in \Pi$ by $(y^\pi_t, u^\pi_t)$. Let $J_T(\pi) = \sum_{t=1}^{T} l_t(y^\pi_t,u^\pi_t), \pi \in \Pi$. Then, the regret for the controller $\mathcal{C}$ is given by
\begin{equation}
R_T(\mathcal{C}) = \mathbb{E}[J_T(\mathcal{C})] - \min_{\pi \in \Pi} \mathbb{E}[J_T(\pi)].
\label{eq:regret-defn}
\end{equation}

\subsection{Assumptions}

We state the assumptions we make below.
\begin{assumption}
The system is stable, i.e., $\rho(A) < 1$. The system matrices $A,B$ are known.
\label{ass:stability}
\end{assumption}
The assumptions on the spectral radius (or the assumption that there is additional knowledge of a feedback rule to stabilize the system) are standard in online learning and control problems \cite{abbasi2011regret, dean2018regret, cohen2019learning, simchowitz2020improper}. We emphasize that analysis without stability or the knowledge of a stabilizing feedback law is still an hard and open challenge in online control. While there are works that investigate simultaneous safe exploration and control such as in Reinforcement Learning \cite{berkenkamp2017safe}, these works do not study the finite performance objective such as regret.  

\begin{assumption}
The noise $w_t$ and $e_t$ are bounded and stochastic i.i.d. Their distribution is known and $\mathbb{E}[w^s_t] = 0$, $\mathbb{E}[e^s_t] = 0$.
\label{ass:noise}
\end{assumption}

\begin{assumption}
The loss function $l_t$ is convex and for $z^\top = [y^\top_t, u^\top_t], (z')^\top = [(y')^\top, (u')^\top]$ such that $R = \max\{\norm{z},\norm{z'},1\}, \norm{l_t(y_t,u_t) -  l_t(y',u')} \leq LR\norm{z - z'}$.
\label{ass:lipschitz}
\end{assumption}
The convexity assumption is standard in online learning and optimization and online control settings. Most of online control especially the setting with general adversarial cost functions and disturbances are built on tools from online convex analysis. This is because the tools for online optimization analysis have been well understood and developed for the convexity setting and such analysis for general non-convex cost setting are still non-existent. The second part of the assumption states that the loss functions are locally Lipschitz. We note that the assumptions stated here are exactly the assumptions in the state-of-the-art work in online control \cite{simchowitz2020improper}. 

\section{Online Control Algorithm}

The online control algorithm we propose for the general controller $\mathcal{C}$ is the {\it adaptive gradient} version of the online DRC (or DRC-GD) proposed in \cite{simchowitz2020improper}. We call this the {\it disturbance response controller - adaptive gradient descent} (DRC-AGD). We briefly review the online DRC in \cite{simchowitz2020improper}, and then present the DRC-AGD algorithm.

\subsection{Online Disturbance Response Controller}

Let's define $y^{nat}$ to be the natural output, the system output when the control inputs are zero, i.e., 
\begin{align}
& y^{nat}_t = e_t + \sum_{s=0}^{t-1} CA^{t-s-1}w_s \nonumber \\
& = y_t - \sum_{s = 1}^{t-1} G^{[s]} u_{t-s}, \ G^{[s]} = CA^{s-1}B. \nonumber
\end{align}

Since $e_t, w_t$ are bounded for all $t$ and $\rho(A) < 1$, $y^{nat}_t$ is bounded for all $t$. We define $R_{nat}$ to be the bound on $y^{nat}_t$. The DRC as defined in \cite{simchowitz2020improper} is parameterized by a $m-$length sequence of matrices, denoted by $M = (M^{[i]})_{i=0}^{m-1}$. The DRC's control decision is given by 
\begin{equation}
u_t = \sum_{s=0}^{m-1} M^{[s]}y^{nat}_{t-s}.
\label{eq:drc}
\end{equation}

Let's define the following class of disturbance response controllers:
\begin{equation}
\mathcal{M}(m,R) = \left\{M = (M^{[s]})_{s=0}^{m-1}: \norm{M} = \sum_s \norm{M^{[s]}} \leq R_M \right\}
\label{eq:drc-class}
\end{equation}
The online learning algorithm or the DRC-GD proposed in \cite{simchowitz2020improper} continuously updates the feedback gain $M$ as the loss functions are revealed. It applies the control input as defined in Eq. \eqref{eq:drc} with the current value of the feedback gain $M$. 
The algorithm then updates the feedback gain $M$ based on the revealed loss function, similar to how the decision is updated in OCO. Thus the disturbance feedback gain $M$ is equivalent to the decision in OCO.

For the choice of regret as defined in Eq. \eqref{eq:regret-defn}, the disturbance response controller is a good choice given that the best disturbance response controller for the realized sequence of cost functions is approximately equal to the best linear dynamic controller. We will show this in the proof of our main result. Thus, by learning the disturbance response controller online the controller can get closer to the optimal linear dynamic controller. We pick the control structure as DRC instead of linear dynamic controller because the DRC control form has advantages from the point of view of online regret analysis. It enables the regret analysis to be approximated by the regret analysis of a limited memory problem, where memory refers to the number of past controller parameters the realized cost at a time $t$ is dependent on. This will not be feasible with the linear dynamic control structure because the control input computed by a linear dynamic controller at any point of time is dependent on the entire history of control inputs unlike Eq. \eqref{eq:drc}.

We introduce the following definitions for ease of presentation. Let $M^{[s]}(j)$ denote the $j$th row of the $M^{[s]}$ matrix. Let $z(i:j)$ denote the sub-vector of the vector $z$ corresponding to the elements from $i$ to $j$. Let $P$ denote the vector given by $P({sq+(j-1)d_y+1:sq+jd_y}) = (M^{[s]}(j))^\top$, where $q = d_yd_u, 1 \leq j \leq d_u$. Essentially, this defines $P$ to be the vector of the transposes of the rows of $M^{[s]}$ stacked one above the other. We introduce the following definitions that will be required for discussing the algorithms.
\begin{definition}
{\it $u_t\left[M_t \vert y^{nat}_{1:t}\right] := \sum_{s = 0}^{m-1} M^{[s]}_t y^{nat}_{t-s}$, \\
$\tilde{y}_t[P_{t:t-h} \vert y^{nat}_{1:t}] := y^{nat}_t + \sum_{s = 1}^h G^{[s]}u_{t-s}$, \\
$F_t\left[P_{t:t-h} \vert y^{nat}_{1:t}\right] := l_t\left(\tilde{y}_t\left[P_{t:t-h} \vert y^{nat}_{1:t}\right], u_t\left[M_t \vert y^{nat}_{1:t}\right]\right)$, \\
$f_t(P \vert y^{nat}_{1:t}) := F_t[\{P,P,...,P\} \vert y^{nat}_{1:t}].$}
\label{def:Ft}
\end{definition}
The term $\tilde{y}_t$ is an approximate output that depends only on the past $h$ control inputs. Consequently this approximate estimate is only a function of $P_{t:t-h}$ for a given $y^{nat}_{1:t}$. The function $F_t$ is the loss $l_t$ evaluated for this approximate output $\tilde{y}_t$ and so it is also only a function of $P_{t:t-h}$. The function $f_t$ is the loss $F_t$ when $P_k$, for all $k$ s.t. $t \geq k \geq t-h$ is fixed to $P$, and so we term it as the memory-less loss. 

Minimizing the regret (Eq. \eqref{eq:regret-defn}) is an Online Convex Optimization problem with Memory (OCO-M) \cite{anava2015online} because the loss function at a time step depends on the past control inputs, which is the case even with the approximated cost $F_t[P_{t:t-h}]$, a function of the truncated output $\tilde{y}_t$. 
Following the key idea in \cite{anava2015online}, the DRC-GD algorithm \cite{simchowitz2020improper} uses the gradient of the memory-less function $f_t(\cdot)$ to update $P$. This, as can be expected, only minimizes the regret of $\sum f_t(\cdot)$ instead of the approximated cost $F_t[P_{t:t-h}]$. But as shown in \cite{anava2015online}, the memory-less regret closely approximates the regret of the approximated cost $F_t[P_{t:t-h}]$, which in turn, as we show later, is a good approximation of the regret of the actual realized cost. 

Let $\mathcal{P}(m,R) = \left\{P : \sum_{s=0}^{m-1} \norm{M^{[s]}} \leq R_M \right\}$. The learning algorithm for the online DRC proposed in \cite{simchowitz2020improper} initializes $P$ to an element drawn from the set $\mathcal{P}(m,R)$. It then updates $P$ along the gradient of the memory-less loss function $f_t(\cdot)$ as the loss functions (or cost) are revealed to continuously improve the feedback controller:
\begin{equation}
P \leftarrow \tn{Proj}_\mathcal{M}\left(P - \eta_{t+1} \partial f_t\left(P \vert y^{nat}_{1:t}\right)\right).
\label{eq:ogd}
\end{equation}

In \cite{simchowitz2020improper}, the authors show that the disturbance response controller with the memory-less gradient update given by Eq. \eqref{eq:ogd}, where $\eta_t$ is fixed to a particular value (see Theorem 2, \cite{simchowitz2020improper}), achieves a regret of $\tilde{\mathcal{O}}(\sqrt{T})$ when the cost functions are general convex functions and polylog$(T)$ when the cost functions are smooth and strongly convex. In this work, we extend this online DRC controller by using an adaptive step rate akin to \cite{hazan2008adaptive} instead of a fixed step rate $\eta$. We discuss our extended algorithm in the next section.

\subsection{Online Disturbance Response Controller: DRC-AGD}

In this section, we present the DRC-AGD algorithm. 
First, we briefly review the adaptive gradient online learning algorithm \cite{hazan2008adaptive} for the standard OCO problem and then present our new regret result for adaptive gradient learning for the OCO-M problem. We then introduce our DRC-AGD online control algorithm and use its result to analyse the regret of the DRC-AGD algorithm.

\subsubsection{Adaptive Gradient Online Learning}

Consider the standard online convex optimization (OCO) setting (see \cite{hazan2008adaptive}). At time $t$, the player chooses an action $u_t$ from some convex subset $\mathcal{K}$ of $\mathbb{R}^n$, where $\max_{x \in \mathcal{K}} \norm{x} \leq D$, and the adversary chooses a convex loss function $f_t(\cdot)$. The regret for the player over duration $T$ is given by 
\begin{equation}
R_T  = \sum_{t=1}^{T} f_t(u_t) - \min_{u \in \mathcal{K}} \sum_{t=1}^{T} f_t(u)
\label{eq:regret-oco}
\end{equation}

Let $f_t$ be $H_t$-strongly convex, i.e., let $f_t(u^{*}) \geq f_t(u) + \nabla f_t (u^{*} - u) + \frac{H_t}{2} \norm{u^{*} - u}^2_2$ and $\norm{\nabla f_t} \leq G_t$. Once the loss function is revealed at time $t$ the algorithm can use the loss function to update its decision. The adaptive gradient online learning algorithm proposed in \cite{hazan2008adaptive} updates the decision $u_t$ by the following gradient step:
\begin{align}
& u_{t+1} = \text{Proj}_\mathcal{K}\left(u_t-\eta_{t+1} \partial\left(f_t(u)+g_t(u)\right)\right) \nonumber \\
& \eta_{t+1} = \frac{1}{\sum H_{1:t}+ \sum \lambda_{1:t}},
\label{eq:aogd}
\end{align}

where $\sum H_{1:t} = \sum_{k=1}^t H_k, \sum \lambda_{1:t} = \sum_{k=1}^t \lambda_k$, and $\lambda_t$s are suitably defined parameters. Here, it is clear that the step rate at each time step is updated by the strong convexity $H_t$ of the loss function at $t$ as defined above. Thus the step rate is adapted and the algorithm is adaptive gradient online learning. The regret for this algorithm can be characterized as in the following Lemma. 
\begin{lemma}
{\it Consider the online update given by Eq. \eqref{eq:aogd} with $g_t(u) = 1/2\lambda_t \norm{u}_2^2$. Then for any sequence of $\lambda_1, \lambda_2,...,\lambda_T$, 
\begin{equation}
R_T \leq \frac{1}{2}D^2 \lambda_{1:T} + \frac{1}{2}\sum_{t=1}^T \frac{(G_t +\lambda_t D)^2}{\sum H_{1:t}+\sum \lambda_{1:t}},
\end{equation}}
\label{lem:aogd-regret}
\end{lemma}

Please see Thoerem 3.1. \cite{hazan2008adaptive} for the proof. This is the basic result that the regret rate results in \cite{hazan2008adaptive} are based on. Here, the parameters $\lambda_{1:T}$ can be suitably chosen based on the convex conditions to achieve intermediate regret rates for intermediate convex conditions of the sequence of loss functions; for example, conditions such as $H_t \propto t^{-\alpha}$. We direct the reader to \cite{hazan2008adaptive} for a more detailed discussion of their results. 

\subsubsection{Adaptive Gradient Online Learning for OCO-M}

In this section we discuss the extension of the adaptive gradient learning to the OCO-M problem. The difference in the OCO-M setting is that the cost function at a particular time $t$ is also dependent on a certain history of the past decisions. More specifically, the cost functions $f_t$ in OCO-M are a function of the decisions upto $h$ time steps in the past, i.e., $u_{t:t-h}$, where $h$ is a given number. Thus, the regret in the OCO-M problem is the following:
\begin{equation}
R_T  = \sum_{t=1}^{T} f_t(u_{t:t-h}) - \min_{u \in \mathcal{K}} \sum_{t=1}^{T} f_t(u),
\label{eq:regret-oco-m}
\end{equation}
where we used $f_t(u)$ as a shorthand notation for the cost when $u_{t-k} = u$, for all $k$, where $0 \leq k \leq h$. In the next theorem we present the equivalent of Lemma \ref{lem:aogd-regret} for the OCO-M problem, which we will use to analyse our main algorithm. 
\begin{theorem}
{\it  For a sequence of $(h+1)$-variate $F_t$ define $f_t(u) = F_t(u,u,...,u)$. Let $G_c$ be an upper bound on the coordinate wise Lipschitz constant of $F_t$, $G_f$ be an upper bound on the Lipschitz constant of $f_t$, $f_t$ be $H_t$-strongly convex, and $D$ be an upper bound on the diameter of $\mathcal{K}$. Consider the online update given by Eq. \eqref{eq:aogd}, with $g_t(u) = 1/2\lambda_t \norm{u}_2^2$. Then for any sequence of $\lambda_1, \lambda_2,...,\lambda_T$, $\lambda_{j} \leq \lambda_{i}, j \geq i$, 
\begin{align}
& R_T = \sum_{t = h+1}^T F_t(u_t,...,u_{t-h}) - \min_{u\in\mathcal{K}} \sum_{t = h+1}^T F_t(u,...,u) \nonumber\\
& \leq \frac{1}{2}D^2 \lambda_{1:T} + \frac{1}{2}\sum_{t=1}^T \frac{\tilde{G}_{f,t}^2}{\sum H_{1:t}+\sum \lambda_{1:t}}, \nonumber
\end{align}
where $\tilde{G}_{f,t} = \sqrt{\left(G_f + \lambda_{t}D\right)(G_f + \lambda_{t}D+2G_ch^{3/2})}$.}
\label{thm:aogd-memory}
\end{theorem}
Please see Appendix for the proof.

\subsubsection{Adaptive Gradient Online Learning for Control}

Here, we extend the adaptive gradient descent learning idea to the online DRC.
The gradient learning algorithm we propose, which we call as DRC-AGD, is the extension of Eq. \eqref{eq:ogd} with an adaptive step rate similar to Eq. \eqref{eq:aogd}:
\begin{align}
& P_{t+1} \nonumber \\
& = \tn{Proj}_\mathcal{P}\left(P_t - \eta_{t+1}\partial \left( \mathbb{E}\left[f_t\left[P_{t} \vert y^{nat}_{1:t}\right]\right] + g_t(P_t)\right)\right) \nonumber \\
& g_t(P) = \frac{1}{2}\lambda_t \norm{P}_2^2, ~ \eta_{t+1} = \frac{1}{\sum H_{1:t}+\sum \lambda_{1:t}},
\label{eq:drc-aogd}    
\end{align}
where the udpate is by the gradient of the memory-less cost $\mathbb{E}\left[f_t\left[P_{t} \vert y^{nat}_{1:t}\right]\right]$, with an adaptive step rate $\eta_{t+1}$, where $H_t$ is the strong convexity of $\mathbb{E}\left[f_t\left[P_{t} \vert y^{nat}_{1:t}\right]\right]$ and $\lambda_t$s are suitably chosen parameters as before.

\begin{algorithm}[H]
\DontPrintSemicolon
\KwInput{Radius $R_M$, and the matrices $G^{[i]}$, $h$.}

Initialize $P_1 \in \mathcal{P}$
 
  \For{t = 1,....,T}    
   { 
	Observe $y_t$ and determine $y^{nat}_t = y_t - \sum_{i = 1}^{t-1} G^{[i]} u_{t-i}$
	
	Choose $u_t = \sum_{s=0}^{m-1} M^{[s]}_t y^{nat}_{t-s}$
	
	Observe the loss function and suffer the loss $l_t(y_t,u_t)$
	
	Set $\eta_{t+1} = \frac{1}{\sum H_{1:t}+\sum \lambda_{1:t}}$
	
	$P_{t+1} = \tn{Proj}_\mathcal{P}\left(P_t - \eta_{t+1}\partial \left( \mathbb{E}\left[f_t\left[P_{t} \vert y^{nat}_{1:t}\right]\right] + \frac{1}{2}\lambda_t \norm{P_t}_2^2\right)\right)$
   }

\caption{Disturbance Response Control - Adaptive Gradient Descent (DRC-AGD)}
\label{alg:drc-agd}
\end{algorithm}
Algorithm \ref{alg:drc-agd} presents the full DRC-AGD algorithm.

\subsubsection{Main Results} In DRC-AGD, the gradient of the memory-less cost $\mathbb{E}\left[f_t\left[P_{t} \vert y^{nat}_{1:t}\right]\right]$ is used. Hence, to apply Theorem \ref{thm:aogd-memory} to the analysis of the DRC-AGD algorithm, we need to establish the strong convexity of $\mathbb{E}\left[f_t\left[P_{t} \vert y^{nat}_{1:t}\right]\right]$. We also need to establish that $G_c$ and $G_f$ exist for the memory-less cost $\mathbb{E}\left[f_t\left[P_{t} \vert y^{nat}_{1:t}\right]\right]$; we prove all of this as part of the main theorem. In the next lemma we characterize the strong convexity of $\mathbb{E}\left[f_t\left[P_{t} \vert y^{nat}_{1:t}\right]\right]$ in terms of the strong convexity $H^l_t$ of $l_t$ (recall how $f_t$ is dependent on $l_t$ in Definition \ref{def:Ft}).

\begin{lemma}
{\it The function $\mathbb{E}\left[f_t\left[P_{t} \vert y^{nat}_{1:t}\right]\right]$ is $H_t$-strongly convex, where 
\begin{equation}
H_t = H^l_t \left(\sigma^2_e + \sigma^2_w \left( \frac{\sigma_{\text{min}}(C)}{1+\norm{A}_2^2}\right)^2\right), \nonumber 
\end{equation} 
$\nabla^2 l_t \geq H^l_t$, $\mathbb{E}[w^{s}_tw^s_t] \geq \sigma^2_w$, $\mathbb{E}[e^{s}_te^s_t] \geq \sigma^2_e$.
}
\label{lem:strongconvexity-F}
\end{lemma}

Please see Proposition 7.1, \cite{simchowitz2020improper} for the proof. 
We introduce an additional definition before we discuss our main theorem. 
\begin{definition}
$\psi(i) = \sum_{j \geq i} \norm{CA^{j-1}B}_2, i > 0$. Since $\rho(A) < 1$, there exists $c > 0$ and $\rho \in (0,1)$ such that $\psi(i) \leq C\rho^i$. $R_{G^{*}} = 1 + \psi(1)$.
\end{definition}

In the next theorem we use Theorem \ref{thm:aogd-memory} to characterize the regret for the DRC-AGD online control algorithm.
\begin{theorem}
{\it Suppose Assumptions \ref{ass:stability}, \ref{ass:noise}, \ref{ass:lipschitz} hold. Suppose the algorithm \ref{alg:drc-agd} is run with $m, h \geq 1$ such that $\psi(m) \leq R_{G^{*}}/T, \psi(h) \leq R_{M}/T$ then
\begin{align}
& R_T(\mathcal{C}) \leq R^2_MR^2_{G^{*}}R^2_{nat}(6L + 4(m+h)) \nonumber \\
& + \frac{1}{2}D^2 \lambda_{1:T} + \frac{1}{2}\sum_{t=1}^T \frac{(\tilde{G}_{f,t})^2}{\sum H_{1:t}+\sum \lambda_{1:t}}, \ \text{where} \nonumber\\
& G_f = G_C = L\sqrt{m}R_MR_{G^{*}}R^2_{nat}, \ D = 2\sqrt{\min\{d_u,d_y\}}R_M, \nonumber
\end{align}
$\tilde{G}_{f,t} = \sqrt{\left(G_f + \lambda_{t}D\right)(G_f+\lambda_{t}D+2G_ch^{3/2})}$. }
\label{thm:drc-agd}
\end{theorem}

Please see the Appendix for the proof. The proof proceeds by splitting the regret (Eq. \eqref{eq:regret-defn}) into several terms; the burn-in loss, algorithm truncation error, f-policy error, comparator truncation error and the policy approximation error. This splitting follows the proof technique in \cite{simchowitz2020improper}. The burn-in loss is just the realized cost corresponding to the first $m+h$ time steps. The burn-in loss can be trivially bounded (see for example Lemma 5.2. \cite{simchowitz2020improper}). The algorithm truncation error is the difference between the realized cost for the remaining horizon and the cost that would be realized with the truncated output approximation $\tilde{y}_t$, i.e., $\sum F_t$. We recall that the output is truncated so that it depends only on the past $h$ control inputs; see Definition \ref{def:Ft} for the truncated output $\tilde{y}_t$ and the corresponding loss $F_t$. This splitting is done because Theorem \ref{thm:aogd-memory} can only be applied to fixed length memory while the actual realized cost is dependent on the entire history of control inputs. The f-policy error is the difference between the cost $\sum F_t$, which is the approximate cost by truncating the memory, and the same cost when $P_k = P ~ \forall ~ k$. Thus, the f-policy error is given by $\sum_{t=m+h+1}^T \mathbb{E}[F_t(P_{t:t-h} \vert y^{nat}_{1:t})] - \inf_{P} \sum_{t=m+h+1}^T \mathbb{E}[f_t(P \vert y^{nat}_{1:t})]$. Given the form of this regret term, we can apply Theorem \ref{thm:aogd-memory} to bound the f-policy error.
We note that the approximated cost with truncated memory under fixed $P$ is different from the realized cost under a fixed disturbance response controller $P$. This introduces the comparator truncation error, the difference of the two costs, i.e., $\inf_{P} \sum_{t=m+h+1}^T \mathbb{E}[f_t(P \vert y^{nat}_{1:t})] - \inf_{P} \sum_{t=m+h+1}^T \mathbb{E}[l_t(y^{P}_t, u^{P}_t)]$. The policy approximation error is the difference between the realized cost for the best fixed $P$ disturbance response controller and the cost for the best linear dynamic controller. The truncation errors and policy approximation error can also be bounded (see \cite{simchowitz2020improper}). We give details of bounding the burn-in loss, truncation errors and the policy approximation error in the Appendix. Putting together the bounds of all these terms gives us the final result.

We note that the regret bound for DRC-AGD has terms similar to the regular adaptive gradient algorithm (see Lemma \ref{lem:aogd-regret}). Given this result, we can apply the analysis similar to \cite{hazan2008adaptive} to establish regret scaling for various convex conditions. In the next corollary we discuss the specific scaling of the regret w.r.t $T$ under various convex conditions and in particular show that the DRC-AGD algorithm interpolates between $T^{1/2}$ and $\log{T}$.

\begin{corollary}
{\it Suppose Assumptions \ref{ass:stability}, \ref{ass:noise}, \ref{ass:lipschitz} hold. Suppose the algorithm \ref{alg:drc-agd} is run with $m, h \geq 1$ such that $\psi(m) \leq R_{G^{*}}/T, \psi(h) \leq R_{M}/T, T \geq 4$ then
\begin{enumerate}
    \item for any sequence of convex loss functions $l_t$
    \begin{equation}
        R_T \leq \tilde{\mathcal{O}}(\sqrt{T}) \nonumber
    \end{equation}
    \item for any sequence of convex loss functions $l_t$ with $H^l_t \geq H$  
    \begin{equation}
        R_T \leq \tilde{\mathcal{O}}(\log{T}) \nonumber
    \end{equation}
    \item for $H^l_t = t^{-\alpha}$, and $0 < \alpha \leq 1/2$
     \begin{equation}
        R_T \leq \tilde{\mathcal{O}}(T^\alpha) \nonumber
    \end{equation} 
     \item for $H^l_t = t^{-\alpha}$, and $\alpha > 1/2$
     \begin{equation}
        R_T \leq \tilde{\mathcal{O}}(\sqrt{T}) \nonumber
    \end{equation} 
\end{enumerate} }
\label{cor:drc-agd}
\end{corollary}

Please see the Appendix for the proof. We see that the DRC-AGD algorithm recovers the $\mathcal{O}(\sqrt{T})$ and $\mathcal{O}(\log{T})$ result for strongly convex and general convex cost functions and at the same time achieves intermediate regret scaling for intermediate convex conditions. We emphasize that the regret scaling of $\tilde{\mathcal{O}}(T^\alpha)$ is valid for a more general condition such as $\sum H_{1:t} \geq t^{1-\alpha}$. 

\section{Conclusion}

In this work we considered the online control of a known linear dynamic system with adversarial disturbances and adversarial cost functions. Our objective is to improve regret rates established for this setting by prior works, which only considered either convex costs or strongly convex costs. Specifically, we addressed the question whether the regret rates can be improved when the convexity of controller cost functions are intermediate, i.e., between strongly convex and convex. 

We proposed an adaptive gradient extension of the disturbance response controller proposed in prior works for the same problem we study. We proved that the proposed online learning controller recovers the previously established regret guarantee of $\mathcal{O}(\sqrt{T})$ for general convex controller cost functions and $\mathcal{O}(\log{T})$ for strongly-convex and smooth controller cost functions (see \cite{simchowitz2020improper}), and achieves an intermediate regret between $\mathcal{O}(\sqrt{T})$ and $\mathcal{O}(\log{T})$ for intermediate convex conditions for the controller cost functions.

\bibliographystyle{plain}
\bibliography{Refs}

\section*{Appendix A: Proof of Theorem \ref{thm:aogd-memory}}
\noindent The regret can be split as 
\begin{align}
& R_T = \sum_{t = h+1}^T F_t(u_t,...,u_{t-h}) - \sum_{t = h+1}^T F_t(u_t,...,u_t) \nonumber \\
& + \sum_{t = h+1}^T F_t(u_t,...,u_t) - \min_{u\in\mathcal{K}} \sum_{t = h+1}^T  F_t(u,...,u) \nonumber \\
& = \sum_{t = h+1}^T F_t(u_t,...,u_{t-h}) - \sum_{t = h+1}^T F_t(u_t,...,u_t) \nonumber \\
& + \sum_{t = h+1}^T f_t(u_t) - \min_{u\in\mathcal{K}} \sum_{t = h+1}^T  f_t(u). \nonumber
\end{align}
Lets call the second term as $\tilde{R}_T$, i.e.,
\begin{equation}
\tilde{R}_T = \sum_{t = h+1}^T f_t(u_t) - \min_{u\in\mathcal{K}} \sum_{t = h+1}^T  f_t(u). \nonumber
\end{equation}
Given that 
\begin{equation}
u_{t+1} = \text{Proj}_\mathcal{K}\left(u_t-\eta_{t+1} \partial\left(f_t(u_t)+g_t(u_t)\right)\right), 
\label{eq:grad-up-xt}
\end{equation}
Lemma \ref{lem:aogd-regret} is applicable to $\tilde{R}_t$. Hence, we have that
\begin{equation}
\tilde{R}_T \leq \frac{1}{2}D^2 \sum \lambda_{1:T} + \frac{1}{2}\sum_{t=1}^T \frac{(G_f +\lambda_t D)^2}{\sum H_{1:t}+\sum \lambda_{1:t}}. \nonumber
\end{equation}  
Next, we bound the first term. By the definition of $G_c$ we have that
\begin{align}
 & \norm{F_t(u_t,...,u_{t-h}) - F_t(u_t,...,u_t)}^2_2 \nonumber \\
 & \leq G^2_c\norm{[u^\top_t,...,u^\top_{t-h}]^\top - [u^\top_t,...,u^\top_{t}]^\top}^2_2 \nonumber \\
 & = G^2_c \sum_{i=1}^{h} \norm{u_t -u_{t-i}}^2_2 \nonumber \\
 & \leq G^2_c \sum_{i=1}^{h} \left(\sum_{j=1}^{i}\norm{u_{t-j+1} -u_{t-j}}_2\right)^2. \nonumber
\end{align}
Using Eq. \eqref{eq:grad-up-xt} we have that
\begin{equation}
\norm{u_{t-j+1} -u_{t-j}}_2 \leq \norm{\eta_{t-j+1} \partial\left(f_{t-j}(u_{t-j})+g_{t-j}(u_{t-j})\right)}_2. \nonumber
\end{equation}
Given that $\norm{\nabla f_{t-j}} \leq G_c$ (this follows from the fact that $L$ is a Lipschitz constant of $f$ iff $\norm{\nabla f}_2 \leq L$ for differentiable $f$), and $\norm{\nabla g_{t-j}(.)}_2 \leq \lambda_{t-j}D$, we have that
\begin{equation}
\norm{u_{t-j+1} -u_{t-j}}_2 \leq \eta_{t-j+1}(G_f + \lambda_{t-j}D). \nonumber
\end{equation}
Using this observation we have that
\[ \norm{F_t(u_t,...,u_{t-h}) - F_t(u_t,...,u_t)}^2_2 \]
\[\leq G^2_c \sum_{i=1}^{h} \left(\sum_{j=1}^{i}\eta_{t-j+1}(G_f + \lambda_{t-j}D)\right)^2 \]
\[\leq G^2_c \sum_{i=1}^{h} \left(\sum_{j=1}^{i}\eta_{t-h}(G_f + \lambda_{t-h}D)\right)^2\]
\[ \leq G^2_c h^3 \eta^2_{t-h} \left(G_f + \lambda_{t-h}D\right)^2. \]
That is 
\begin{align}
& \norm{F_t(u_t,...,u_{t-h}) - F_t(u_t,...,u_t)}_2 \nonumber \\
& \leq G_ch^{3/2}\eta_{t-h}\left(G_f + \lambda_{t-h}D\right). \nonumber 
\end{align}
Hence 
\begin{align}
& \sum_{t = h+1}^T F_t(u_t,...,u_{t-h}) - \sum_{t = h+1}^T F_t(u_t,...,u_t) \nonumber \\
& \leq G_ch^{3/2} \sum_{t = h+1}^T \eta_{t-h}\left(G_f + \lambda_{t-h}D\right) \nonumber \\
& = G_ch^{3/2} \sum_{t = h+1}^T \frac{\left(G_f + \lambda_{t-h}D\right)}{\sum H_{1:t-h}+\sum \lambda_{1:t-h}} \nonumber \\
& \leq G_ch^{3/2} \sum_{t = 1}^T \frac{\left(G_f + \lambda_{t}D\right)}{\sum H_{1:t}+\sum \lambda_{1:t}}. \nonumber 
\end{align}
Combining this with the bound on $\tilde{R}_T$ we get that
\begin{equation}
R_T \leq \frac{1}{2}D^2 \sum \lambda_{1:T} + \frac{1}{2}\sum_{t=1}^T \frac{\tilde{G}_{f,t}^2}{\sum H_{1:t}+\sum \lambda_{1:t}}, \nonumber
\end{equation}
where $\tilde{G}_{f,t} = \sqrt{\left(G_f + \lambda_{t}D\right)(G_f + \lambda_{t}D+2G_ch^{3/2})}$. $\blacksquare$



\section*{Appendix C: Proof of Theorem \ref{thm:drc-agd}}

\noindent 
For a policy $\pi^{*} \in \Pi$
\begin{equation}
J_T(\mathcal{C}) - J_T(\pi^{*}) = \sum_{t=1}^T l_t(y_t,u_t) - \sum_{t=1}^T l_t(y^{\pi^{*}}_t, u^{\pi^{*}}_t) \nonumber 
\end{equation} 
We can split the regret as in \cite{simchowitz2020improper}:
\[\mathbb{E}[J_T(\mathcal{C})] - \mathbb{E}[J_T(\pi^{*})] = \underbrace{\sum_{t=1}^{m+h} \mathbb{E}[l_t(y_t,u_t)]}_{\text{burn-in loss}} \]
\[+ \underbrace{\sum_{t=m+h+1}^T \mathbb{E}[l_t(y_t,u_t)] - \sum_{t=m+h+1}^T \mathbb{E}[F_t(P_{t:t-h} \vert y^{nat}_{1:t})]}_{\text{algorithm truncation error}}\]
\[+ \underbrace{\sum_{t=m+h+1}^T \mathbb{E}[F_t(P_{t:t-h} \vert y^{nat}_{1:t})] - \inf_{P} \sum_{t=m+h+1}^T \mathbb{E}[f_t(P \vert y^{nat}_{1:t})]}_{\text{f-policy error}} \]
\[ + \underbrace{\inf_{P} \sum_{t=m+h+1}^T \mathbb{E}[f_t(P \vert y^{nat}_{1:t})] - \inf_{P} \sum_{t=m+h+1}^T \mathbb{E}[l_t(y^{P}_t, u^{P}_t)]}_{\text{comparator truncation error}} \]
\[ + \underbrace{\inf_{P} \sum_{t=1}^T \mathbb{E}[l_t(y^{P}_t, u^{P}_t)] - \sum_{t=1}^T \mathbb{E}[l_t(y^{\pi^{*}}_t, u^{\pi^{*}}_t)].}_{\text{policy approximation error}}. \]

We leverage the results from \cite{simchowitz2020improper} to bound the following terms: (i) burn-in loss, (ii) algorithm truncation error, (iii) comparator truncation error, and (iv) policy approximation error. From Lemma 5.2, \cite{simchowitz2020improper}, we have that 
\begin{equation}
\sum_{t=1}^{m+h} \mathbb{E}[l_t(y_t,u_t)] \leq 4R^2_{G^{*}}R^2_{nat}R^2_M(m+h). \nonumber
\end{equation}
From Lemma 5.3, \cite{simchowitz2020improper}, we have that
\begin{equation}
\mathbb{E}[\text{Truncation errors}] \leq 4LTR_{G^{*}}R^2_{nat}R^2_M\psi(h+1) \nonumber 
\end{equation}
Finally, from Theorem 1, \cite{simchowitz2020improper}, we have that
\begin{equation}
\mathbb{E}[\text{Policy app. error}] \leq 2LTR_MR^2_{G^{*}}R^2_{nat}\psi(m) \nonumber
\end{equation}
Next we bound the f-policy error term. Theorem \ref{thm:aogd-memory} applies to this term. From Lemma 5.4, \cite{simchowitz2020improper}, we have that $f_t(.\vert y^{nat}_t)$ is $G_f$-Lipschitz, where $G_f = L\sqrt{m}R_MR_{G^{*}}R^2_{nat}$, $F_t(.\vert y^{nat}_t)$ is $G_f$-Lipschitz coordinate wise, i.e., $G_c = G_f$, and $D = 2\sqrt{\min\{d_u,d_y\}}R_M$. Then applying Theorem \ref{thm:aogd-memory} to the f-policy error term we get that
\begin{equation}
\mathbb{E}[\text{f-policy error}] \leq \frac{1}{2}D^2 \sum \lambda_{1:T} + \frac{1}{2}\sum_{t=1}^T \frac{(\tilde{G}_{f,t})^2}{\sum H_{1:t}+\sum \lambda_{1:t}}. \nonumber
\end{equation}
This completes the proof.
$\blacksquare$

\section*{Appendix D: Proof of Corollary \ref{cor:drc-agd}}

\noindent Consider the term 
\begin{equation}
\hat{R}_T = \frac{1}{2}D^2 \sum \lambda_{1:T} + \frac{1}{2}\sum_{t=1}^T \frac{(\tilde{G}_{f,t})^2}{\sum H_{1:t}+\sum \lambda_{1:t}}. \nonumber
\end{equation}
We make the following observation.
\begin{align}
& (\tilde{G}_{f,t})^2 = \left(G_f + \lambda_{t}D\right)\left(G_f+\lambda_{t}D+2G_ch^{3/2}\right) \nonumber \\
& \leq 2 \left(G_f + \lambda_{t}D\right)^2 + 2 G^2_ch^{3} \leq 4G^2_f + 4\lambda^2_{t}D^2 + 2 G^2_ch^{3}. \nonumber
\end{align}
Hence,
\begin{align}
& \hat{R}_T \leq \frac{1}{2}D^2 \sum \lambda_{1:T} + \frac{1}{2}\sum_{t=1}^T \frac{4G^2_f + 4\lambda^2_{t}D^2 + 2 G^2_ch^{3}}{\sum H_{1:t}+\sum \lambda_{1:t}} \nonumber \\
& \leq \frac{1}{2}D^2 \sum \lambda_{1:T} + 2\sum_{t = 1}^T\lambda_{t}D^2 + \sum_{t=1}^T \frac{2G^2_f + G^2_ch^{3}}{\sum H_{1:t}+\sum \lambda_{1:t}} \nonumber \\
& \leq \frac{5}{2}D^2 \sum \lambda_{1:T} + \sum_{t=1}^T \frac{2G^2_f + G^2_ch^{3}}{\sum H_{1:t}+\sum \lambda_{1:t}}. 
\label{eq:hR-bound}
\end{align}
Let $\hat{G}^2_f := 2G^2_f + G^2_ch^{3}$. Next, we prove the main results case by case.

\noindent {\em Case 1, any sequence of convex $l_t$}: For this case set $\lambda_1 = \sqrt{T}$ and $\lambda_t = 0$, $t \geq 2$. Then from Theorem \ref{thm:drc-agd} and Eq. \eqref{eq:hR-bound} we get that 
\begin{align}
& R_T(\mathcal{C}) \leq R^2_MR^2_{G^{*}}R^2_{nat}(6L + 4(m+h)) \nonumber \\
& + \frac{5}{2}D^2 \sum \lambda_{1:T} + \sum_{t=1}^T \frac{\hat{G}_{f}^2}{\sum H_{1:t}+\sum \lambda_{1:t}} \nonumber \\
& \leq R^2_MR^2_{G^{*}}R^2_{nat}(6L + 4(m+h)) + \frac{5}{2}D^2\sqrt{T} + \nonumber \\
& + \hat{G}_{f}^2\sum_{t=1}^T \frac{1}{\sqrt{T}} = \mathcal{O}(\sqrt{T}). \nonumber
\end{align}

\noindent {\em Case 2, any sequence of convex $l_t$ such that $H^l_t \geq H$}: In this case, from Lemma \ref{lem:strongconvexity-F}
\begin{align}
H_t \geq H^l_t\left(\sigma^2_e + \sigma^2_w \left( \frac{\sigma_{\text{min}}(C)}{1+\norm{A}_2^2}\right)^2\right) \nonumber \\
\geq H\left(\sigma^2_e + \sigma^2_w \left( \frac{\sigma_{\text{min}}(C)}{1+\norm{A}_2^2}\right)^2\right) = \tilde{H}. \nonumber
\end{align}
Set $\lambda_t = 0$, then from Theorem \ref{thm:drc-agd} and Eq. \eqref{eq:hR-bound} we get that 
\begin{align}
& R_T(\mathcal{C}) \leq R^2_MR^2_{G^{*}}R^2_{nat}(6L + 4(m+h)) \nonumber \\
& + \frac{5}{2}D^2 \sum \lambda_{1:T} + \sum_{t=1}^T \frac{\hat{G}_{f}^2}{\sum H_{1:t}+\sum \lambda_{1:t}} \nonumber \\
& \leq R^2_MR^2_{G^{*}}R^2_{nat}(6L + 4(m+h)) + \hat{G}_{f}^2\sum_{t=1}^T \frac{1}{t\tilde{H}} \nonumber \\
& = \mathcal{O}(\log{T}). \nonumber 
\end{align}
\noindent {\em Case 3, $H^l_t = Ht^{-\alpha}$, and $0 < \alpha \leq 1/2$}: From Lemma \ref{lem:strongconvexity-F} $H_t \geq \tilde{H}t^{-\alpha}$. Set $\lambda_1 = \tilde{H}T^{\alpha}$ and $\lambda_t = 0$, $t > 1$. Then from Theorem \ref{thm:drc-agd} and Eq. \eqref{eq:hR-bound} we get that 
\begin{align}
& R_T(\mathcal{C}) \leq R^2_MR^2_{G^{*}}R^2_{nat}(6L + 4(m+h)) \nonumber \\
& + \frac{5\tilde{H}}{2}D^2T^{\alpha} + \frac{\hat{G}_{f}^2}{\tilde{H}}\sum_{t=1}^T \frac{1}{\left(\sum_{k = 1}^tk^{-\alpha} + T^{\alpha}\right)} \nonumber.
\end{align}
Now $\sum_{k = 1}^tk^{-\alpha} \geq \int_{0}^{t-1} (u+1)^{-\alpha} du = (1-\alpha)^{-1}(t^{1-\alpha}-1)$. Using this fact we get that
\begin{align}
& R_T(\mathcal{C}) \leq R^2_MR^2_{G^{*}}R^2_{nat}(6L + 4(m+h)) \nonumber \\
& + \frac{5\tilde{H}}{2}D^2T^{\alpha} + \frac{\hat{G}_{f}^2(1-\alpha)}{\tilde{H}}\sum_{t=1}^T t^{\alpha-1} \nonumber \\
& \leq R^2_MR^2_{G^{*}}R^2_{nat}(6L + 4(m+h)) \nonumber \\
& + \frac{5\tilde{H}}{2}D^2T^{\alpha} + \frac{\hat{G}_{f}^2(1-\alpha)}{\tilde{H}\alpha}T^{\alpha} = \mathcal{O}(T^{\alpha}). \nonumber
\end{align}
\noindent {\em Case 4, $H^l_t = Ht^{-\alpha}$, and $\alpha \geq 1/2$}: In this case too $H_t \geq \tilde{H}t^{-\alpha}$. Set $\lambda_1 = \tilde{H}T^{1/2}$ and $\lambda_t = 0$, $t > 1$. Then from Theorem \ref{thm:drc-agd} and Eq. \eqref{eq:hR-bound} we get that 
\begin{align}
& R_T(\mathcal{C}) \leq R^2_MR^2_{G^{*}}R^2_{nat}(6L + 4(m+h)) \nonumber \\
& + \frac{5\tilde{H}}{2}D^2T^{1/2} + \frac{\hat{G}_{f}^2}{\tilde{H}}\sum_{t=1}^T \frac{1}{T^{1/2}} = \mathcal{O}(\sqrt{T}) \blacksquare \nonumber
\end{align}

\end{document}